\documentclass{amsart}

\usepackage{etex}
\usepackage{amsmath, amssymb}
\usepackage{array}
\usepackage[frame,cmtip,arrow,matrix,line,graph,curve]{xy}
\usepackage{graphpap, color, paralist, pstricks}
\usepackage[mathscr]{eucal}
\usepackage[pdftex]{graphicx}
\usepackage[pdftex,colorlinks,backref=page,citecolor=blue]{hyperref}
\usepackage{pifont}
\usepackage{cancel}
\usepackage{mathtools}
\usepackage{tikz}

\usepackage{bbm}
\usepackage{dsfont}

\newcommand{\cF}{\mathcal{F} }

\newcommand{\cI}{\mathcal{I} }

\newcommand{\pr}{\mathbb{P}}

\newcommand{\E}[0]{\mathbb{E}}

\newcommand{\beq}[1]{\begin{equation}\label{#1}}
\newcommand{\enq}[0]{\end{equation}}

\newcommand{\mn}[0]{\medskip\noindent}
\newcommand{\nin}[0]{\noindent}

\newcommand{\sub}[0]{\subseteq}
\newcommand{\sm}[0]{\setminus}
\renewcommand{\dots}[0]{,\ldots,}

\newcommand{\un}[0]{\underline}

\newcommand{\cee}[0]{{\mathcal C}}
\newcommand{\D}[0]{{\mathcal D}}

\newcommand{\f}[0]{{\mathcal F}}

\newcommand{\I}[0]{{\mathcal I}}

\newcommand{\T}[0]{{\mathcal T}}

\newcommand{\Y}[0]{{\mathcal Y}}

\newcommand{\Z}[0]{{\mathcal Z}}

\newcommand{\Ra}[0]{\Rightarrow}

\newcommand{\ww}{\mbox{{\sf w}}}

\newcommand{\0}[0]{\emptyset}

\newcommand{\C}[2]{{{#1}\choose{{#2}}}}

\newcommand{\ga}[0]{\alpha }
\newcommand{\gb}[0]{\beta }
\newcommand{\gc}[0]{\gamma }
\newcommand{\gd}[0]{\delta }
\newcommand{\gD}[0]{\Delta }
\newcommand{\gG}[0]{\Gamma }

\newcommand{\gl}[0]{\lambda }

\newcommand{\go}[0]{\omega}
\newcommand{\gO}[0]{\Omega}

\newcommand{\gs}[0]{\sigma}

\newcommand{\gz}[0]{\zeta}
\newcommand{\eps}[0]{\varepsilon }
\newcommand{\vt}[0]{\vartheta}

\newcommand{\vr}[0]{\varrho}
\newcommand{\vp}[0]{\varphi}

\newcommand{\pt}{\vp}

\newcommand{\prh}[1][]{\pr_h}

\newtheorem{thm}{Theorem}[section]
\newtheorem{prop}[thm]{Proposition}

\newtheorem{lemma}[thm]{Lemma}
\newtheorem{obs}[thm]{Observation}
\newtheorem{cor}[thm]{Corollary}

\theoremstyle{definition}

\newtheorem{claim}[thm]{Claim}

\newtheorem{question}[thm]{Question}

\begin{document}

\title{Linear cover time is exponentially unlikely}

\author{Quentin Dubroff \and Jeff Kahn}
\thanks{Department of Mathematics, Rutgers University}
\thanks{JK was supported by NSF Grant DMS1954035}
\email{qcd2@math.rutgers.edu,jkahn@math.rutgers.edu}
\address{Department of Mathematics, Rutgers University \\
Hill Center for the Mathematical Sciences \\
110 Frelinghuysen Rd.\\
Piscataway, NJ 08854-8019, USA}

\begin{abstract}

Proving a 2009 conjecture of Itai Benjamini, we show:

\mn
\textbf{Theorem} 
\emph{For any C there is an $\eps>0$ such that for any simple graph G on $V$ of size $n$,
and $X_0,\ldots$ an ordinary random walk on G,}
\[
\pr(\{X_0\dots X_{Cn}\}= V) < e^{-\eps n}.
\]

\mn
A first ingredient in the proof of this is a similar statement for Markov 
chains in which all transition probabilities are sufficiently small relative to $C$.

\end{abstract}

\maketitle


\section{Introduction}\label{Intro}

We are motivated by a surprisingly basic 
question that we first heard from Ori Gurel--Gurevich
in 2010:  \emph{is it true that for any fixed $C$ and $n$-vertex 
simple graph $G$, the probability that a random walk on $G$ covers $V(G)$ 
in $Cn$ steps is exponentially small in $n$?}  
(Some usage notes are included at the end of this section.)

A positive answer was conjectured by Itai Benjamini in 2009 (\cite{IB}; see also \cite{BGM}),
and given some support by 
a quite amazing argument of Benjamini, Gurel-Gurevich,
and Morris \cite{BGM}, showing that the answer is yes if we assume 
any fixed bound $\gD$ on the maximum
degree of $G $ (with the constant in the exponent then depending on $\gD$
as well as $C$).  That the answer is yes for trees was shown 
by Yehudayoff \cite{Yehudayoff}, who also observed that when $G$ is an expander,
a positive answer follows easily from the (less easy) large deviation bound of Gillman \cite{Gillman}. 

Here we answer the question:
\begin{thm}\label{thmMT}
For any C there is an $\eps>0$ such that for any simple graph G on $V$ of size $n$,
and $X_0,\ldots$ an ordinary random walk on G (with any rule for $X_0$),
\[
\pr(\{X_0\dots X_{Cn}\}= V) < e^{-\eps n}.
\]
\end{thm}

The machine underlying the proof of Theorem~\ref{thmMT} is the following statement for 
general Markov chains, which seems of independent interest.
Here and in the corollary that follows, $(X_t)$ is a Markov chain on $V$
(of size $n$)
with transition probabilities $\pt(\cdot,\cdot)$, and 
$X_I$ is the \emph{set} $\{X_t:t\in I\}$ (for a set of indices $I$).
\begin{thm}\label{main}
For each $C$ and $\gb>0$ there is a $\gd =\gd(C,\gb)>0$ such that
if $W\sub V$, $|W|>\gb n$, $M\leq C|W|$, and
\beq{vpvw}
\max\{\pt(v,w): v\in V, v\neq w\in W\} <\gd,
\enq
then
\[
\pr(X_{[M]} \supseteq W) = \exp[-\gO_{C,\gb}|W|].
\]
\end{thm}

\nin
In particular (roughly) the conclusion of 
Theorem~\ref{thmMT} holds for \emph{any} Markov chain in which the transition
probabilities are small enough relative to $C$.  This includes 
Theorem~\ref{thmMT} itself when the minimum degree of $G$ 
is sufficiently large: 

\begin{cor}\label{CorLD}
For each C there are $\eps>0$ and d such that for any RW $(X_t)$ on a graph of 
minimum degree at least $d$,
\[
\pr(X_{[Cn]}=V) < e^{-\eps n}.
\]
\end{cor}
\nin
This again seems interesting in its own right; e.g., we don't know another way to prove
Theorem~\ref{thmMT} even for the Hamming Cube ($\{0,1\}^m$ with the natural adjacencies),
the scene of some of our early skirmishes with the present problem.

One might hope that Theorem~\ref{thmMT} could now be handled by some combination of 
Corollary~\ref{CorLD} and the ideas of \cite{BGM}, but this seems to be a dead end.
(We did at least manage to ``borrow'' \cite{BGM}'s title.
The only antecedent we know of
for what follows
is \cite[Lemma 2]{Aldous}, whose sibling,
the present Lemma~\ref{assoc}, was our starting point.
In particular, beautiful work of \cite{DLP,Ding,Zhai}, showing (see \cite[Theorem 1.1]{Zhai})
that cover time is ``exponentially concentrated'' in a different sense,
seems
unconnected to what we do here.)

The actual proof of Theorem~\ref{thmMT} is 
based especially on the following easy consequence of Theorem~\ref{main},
which again applies to general Markov chains
(and in which $\vp_W$ refers to the ``induced'' chain on $W$; see ``Usage'' below).

\begin{cor}\label{corP}
With $|V|=n$, suppose the partition
$V^0\cup \cup_{i=1}^k V_i$ of $V$, 
and $U_i\sub V_i$ ($i\in [k]$), satisfy
$|V_i|> \vt n$ $\forall i$; $|V^0| < (1-\gc) n$; $|U_i|> \gc |V_i|$; and,
with $\pt_i=\pt_{_{V_i}}$,
\beq{Ui}
\max\{\pt_i(v,w): v\in V_i, v\neq w\in U_i\} <\gd(C/\gc^2,\gc).
\enq
Then
\beq{implied}
(\pr(X_{[Cn]} =V) \leq )  \,\,\,\, \pr(X_{[Cn]}\supseteq \cup V_i)= \exp[-\gO(n)],
\enq
where the implied constant depends on the constants $C,\gc$ and $\vt$.
\end{cor}

\nin
\emph{Proof.}
Since $|\cup V_i|>\gc n$, if $X_{[Cn]}\supseteq \cup V_i$ then there is an $i$ 
for which the first $(C/\gc)|V_i|$ steps 
of the induced chain on $V_i$ cover $U_i$,
 a set of size at least $\gc |V_i|$.
So Theorem~\ref{main} bounds the l.h.s.\ of \eqref{implied} by
$ke^{-\gO(n)} =e^{-\gO(n)}$.

\qed

In what follows $\gc$ will be a ``true'' constant, meaning one not depending on $C$,
and $\vt$ will be a function of $C$; thus the implied constant in \eqref{implied}
depends only on $C$ and we have Theorem~\ref{thmMT} whenever we can show the existence of 
the desired partition.
Of course not every Markov chain admits such a partition
(or we would have the nonsensical claim that Theorem~\ref{thmMT} holds for general chains), 
but it seems possible that RW (again, on a simple $G$) does.

We will find it convenient to set (for the rest of the paper) 
\[   
\gc = 0.1,
\]   
but stress that any slightly small 
``true'' constant would do as well.

\begin{question}\label{Question}
Is it true that for each $C$ there is a $\vt$ for which, for any $G$, RW on $G$ admits
a partition as in Corollary~\ref{corP}?
\end{question}

\nin
(Of course for Theorem~\ref{thmMT} it would be enough to have a positive answer with
$\gc$ also a function of $C$.)

In the event, we are only able to produce (more accurately, show existence of) such 
a partition under a pair of restrictions on $G$, but can also show that if either of these
is violated then Theorem~\ref{thmMT} holds for other reasons.
Failure of the first restriction, which forbids too many large degrees, is handled by
the next lemma, which may be thought of (not quite accurately because of the difference
in the degree bounds) as a substantial extension of Corollary~\ref{CorLD}.


\begin{lemma}\label{lemmaB}  
For each D there is a $\gD$ such that if
\beq{LBineq}
|\{v: d_G(v)> \gD\}| > \gc n,
\enq
then $\pr(X_{[Dn]}= V) =\exp[-\gO (n)]$.
\end{lemma}

\nin 
We postpone specifying the second restriction, which will be easier to do in the context of 
Section~\ref{Partitions}
(see Lemma~\ref{lemmaY} and \eqref{Dsmall}).

\mn
\emph{Remark.}
Lemma~\ref{lemmaB} is the only place where we use simplicity of $G$, the
rest of what we do being valid for general reversible chains.
At that level, Theorem~\ref{thmMT} does not hold without some restriction, but e.g., 
the argument of Section~\ref{Partitions} goes through essentially 
unchanged to show (with $\pi$ denoting stationary distribution): 
\begin{thm}
Let $(X_i)_{i\geq 0}$ be a reversible Markov chain on $V$, and suppose there exists $W \subseteq V$ with $|W| \geq \ga n$ and $\pi_v\leq K\pi_w~\forall v,w\in W$.  Then
\[\pr(X_{[Dn]} = V) = \exp[-\gO(n)],\]
where the implied constant depends on $\ga$, $K$, and $D$.
\end{thm}

Before closing this discussion we mention an obvious challenge:
\begin{question}\label{Qrate}
Can one say anything reasonable about the value of $\eps$
in Theorem~\ref{thmMT}?
\end{question}
\nin
Whatever value can be extracted from our argument will be quite bad
(We suspect it's not as bad as what could be gotten from \cite{BGM}, but are not volunteering
to make this comparison.)  As far as we know, it could be that, for slightly large $C$, 
complete graphs---for which $\pr(\mbox{cover})$ is roughly $\exp[-e^{-C}n]$---are 
more or less the worst case;
but note that for $C=1$ (e.g.), the probability is larger for a path.
At any rate, given how far we are from a decent value, there's clearly no point
in trying to optimize anything, and we instead do what we can to keep things reasonably simple.

\mn
\emph{Outline.}
Following brief preliminaries in Section~\ref{Preliminaries},
Theorem~\ref{main} and Lemma~\ref{lemmaB} are proved in Sections~\ref{mainsec}
and \ref{largeB} (respectively), and the derivation of Theorem~\ref{thmMT} via Corollary~\ref{corP}
is given in Section~\ref{Partitions}.
To give some sense of
Corollary~\ref{corP}, 
two ``bonus'' sections at the end of the paper return to known
cases of Theorem~\ref{thmMT} for which our machinery operates 
relatively simply:
Section~\ref{expand} treats expanders,
and 
might be read as an interlude following
Observation~\ref{Bvobs}.
Section~\ref{Trees}, which reproves Yehudayoff's result for trees,
can be read at any point (including this point).
Finally, we have added an appendix treating a martingale concentration statement
related to Section~\ref{mainsec} (see following \eqref{H-A}),
which is not needed for present purposes but might be of independent interest.

\mn
\emph{Usage.}
We consider Markov chains $(X_t)_{t\geq 0}$ on default state space $V$ of size $n$,
as usual using $\pi$ for stationary distribution. 
We use $\pr_v$ and $\E_v$ for probability and expectation given $X_0 = v$. 
For $B\subseteq V$, the \emph{hitting time} of $B$ is 
$T_B = \min\{t : X_t\in B\}$, and 
the \emph{positive hitting time} is $T_B^+ := \min\{t>0 : X_t\in B\}$
($=T_B$ if $X_0\not\in B$).

We use $\vp$ for transition probabilities 
and $\vp_{_W}$ for
transition probabilities in the \textit{induced chain} on $W \subseteq V$; that is,
$\pt_W(u,v) = \pr_u(X_{T_W^+} = v)$. (This usage is not universal; 
e.g.\ \cite{AF} uses ``chain watched on $W$" here and
``induced chain'' differently.)

Throughout $G =(V,E)$ is a (finite, connected) \emph{simple} graph, with (again) $|V|=n$.
Usage here is pretty standard:  $N_v$ for the neighborhood of (i.e.\ vertices 
adjacent to) $v$; $d_v$---or, if necessary, $d_G(v)$---for $|N_v|$ (the \emph{degree} of $v$);
and, for $X\sub V$, $N(X) = \cup_{x\in X} N_x$. 
We use \emph{random walk (RW) on} $G$ for a Markov chain on $V$ with 
$\vp(v,w) = \mathbbm{1}_{\{v\sim w\}}/d_v$ (with any choice of $X_0$),
recalling that then $\pi_v = d_v/(2|E|)$.

We use $[n]$ for $\{1,2,\ldots, n\}$ and always assume $n$ is large enough 
to support our arguments. 
To avoid needless clutter, we allow a few irrelevant abuses such as (usually)
pretending large numbers are integers.

\mn

\section{Preliminaries}\label{Preliminaries}

We collect here only a few items that will be needed below (and that most 
readers might profitably skip).  
For general background on both Markov chains and martingales, see e.g.\ 
\cite{MCMT}.

Recall that a Markov chain (with stationary distribution $\pi$) is \textit{reversible} if, for any $u,v \in V$,
\[   
\pi_u \pt(u,v) = \pi_v \pt(v,u);
\]   
equivalently: for any $v_0,\ldots,v_k\in V$,
\beq{pathrev}
\pi_{v_0}\pt(v_0,v_1)\pt(v_1,v_2)\cdots\pt(v_{k-1},v_k) = \pi_{v_k}\pt(v_k,v_{k-1})\pt(v_{k-1},v_{k-2})\cdots\pt(v_{1},v_0).
\enq
(A reversible Markov chain is the same thing as RW on a weighted graph----that is,
with weights $\ww(\cdot,\cdot)$ on edges and $\vp(v,w)\propto \ww(v,w)$---but we 
won't need this.)

The next two inequalities are for use in Section~\ref{Partitions}.
The first bounds transition probabilities in terms of return probabilities.
The second---monotonicity of return probabilities---will be used to 
deal with a \emph{tiny} technical annoyance,

\begin{lemma}\cite[Lemma 3.20]{AF}\label{AF320}
For any two states $v$ and $w$ of a reversible Markov chain (and any s,t),
\[
\frac{\pt^{t+s}(v,w)}{\pi_v} \leq \bigg[\frac{\pt^{2t}(v,v)}{\pi_v}\frac{\pt^{2s}(w,w)}{\pi_w}\bigg]^{1/2}.
\]
\end{lemma}

\begin{lemma}\cite[Proposition 10.25]{MCMT}\label{mono}
For any state $v$ of a reversible Markov chain (and any t),
\[\pt^{2t+2}(v,v) \leq \pt^{2t}(v,v).\]
\end{lemma}

The following basic martingale facts will be used in the proof of Theorem~\ref{main}
(in Section~\ref{mainsec}).
The first is a weak form of the Martingale Convergence Theorem; see e.g. \cite[Theorem~5.1]{KT} and \cite[Proposition~A.11(i)]{MCMT}. 

\begin{thm}\label{conv}
If $X_s$ is a supermartingale with (for some L) $|X_s|\leq L$ for all $s\geq 0$, then 
there is a random variable $X$ such that $X_s$ converges to $X$ with probability one, and
\beq{EXn}
\E X_s \leq \E X \,\,\,\forall s.
\enq
\end{thm}

\nin 
All limits in Section~\ref{mainsec} are easily seen to exist everywhere, so for us
the important part of Theorem~\ref{conv} 
is \eqref{EXn}.

Lastly, we recall (a special case of) the ``Hoeffding-Azuma'' Inequality:  
\begin{thm}\cite[Theorem A.10]{MCMT}\label{HA}
If $X_s$ is a martingale with $|X_{s+1}- X_s| \leq L$ for all $s\geq 0$, then
\[\pr(X_k - \E X_k > \eta) \leq e^{-\eta^2/(2kL^2)}.\]
\end{thm}

\section{Proof of Theorem \ref{main}}\label{mainsec}

As mentioned above, our initial inspiration was Aldous' paper \cite{Aldous}.  Our notation here
is more or less his, and Lemma~\ref{assoc} was suggested by his Lemma~2.

The proof of Theorem~\ref{main}, given at the end of this section, is a simple 
application of the material we are about to develop.  
Until then we keep the discussion slightly more general---if not as general as it 
might have been---to support a second application in the proof of Lemma~\ref{lemmaB} in Section~\ref{largeB}.

We consider some $W\sub V$ and hope to show, under suitable assumptions,
that 
\beq{coverW}
\pr(X_{[M]}\supseteq W) < e^{-\gO(n)},
\enq
where the implied constant depends on $C$ and $\gb <|W|/n$.

We will work with a parameter $K$, a (slightly large) function of $C$ and $\gb$;  
but as the value of $K$ plays no role in the first half (or so) of this discussion, 
we leave it unspecified until---in Lemma~\ref{union}---it becomes relevant.

Given $W\subseteq V$, let $|W|=m$ and 
\beq{lambda}
\gl = 1- \max\{\pt(v,w):v\in V, v\neq w\in W\} .
\enq
(Though we've kept track of $\gl$ here, in our applications it will be at least 
$1/2$ and its precise value will be unimportant.)

Let $L = e^K$ and define random variables
\[
H_v(t) = \prod_{i=0}^t [1 - \pt(X_i, v)]
\]
(so $H_v(-1)=1$) and
\[
r_v = \min\{t : H_v(t) < \gl/L\}.
\]

\nin
Write $a\wedge b$ for $\min(a,b)$ and parse $a\wedge b-1 = (a\wedge b)-1$.
Define martingales 
\[    
\xi_s^v = \mathbbm{1}_{\{T_v > s \wedge r_v\}}/ H_v( s \wedge r_v - 1)
\]   
and
\[   
\xi_s=\xi_s^W = \sum_{w \in W} \xi_s^w.
\]   
We omit the (standard, easy) proof that they \emph{are} martingales
(see e.g.\ the proof of \cite[Lemma 1]{Aldous} 
for essentially the same argument).  We assume (as we may) that $X_0\not\in W$, so 
\[
\E\xi_M=\xi_0 = m,
\]
and observe that
\beq{xisxis}
|\xi_s - \xi_{s-1}| \leq L/\gl^2.
\enq
[\emph{Because}:  with sums over $v$'s with $ T_v, r_v > s-1$ (i.e.\ those that can contribute here), we have

\nin
\begin{align*}
    |\xi_{s} - \xi_{s-1}| \,\,
    &\mbox{$\,\,= \,\,|\sum (\xi_s^v - \xi_{s-1}^v)| \,\,
    \leq \,\,\max\{L/\lambda, \sum \xi_{s-1}^v/(1-\pt(X_{s-1},v)) - \xi_{s-1}^v\}$}\\
    &\mbox{$\,\,= \,\,\max\{L/\lambda, \sum\frac{\pt(X_{s-1},v)\xi_{s-1}^v}{1-\pt(X_{s-1},v)}\}\,\,
    \leq \,\,\max\{ L/\lambda,  (1/\gl)\sum\pt(X_{s-1},v)\xi_{s-1}^v\}$}\\
    &\mbox{$\,\leq \,\,\max\{ L/\lambda, (L/\gl^2) \sum\pt(X_{s-1},v)\}\,\,
    \leq L/\gl^2$.]}
\end{align*}

The Hoeffding-Azuma Inequality, 
Theorem \ref{HA}, thus gives
\beq{H-A}     
\pr(\xi_M  < m/2)  = e^{-\gO(m)}.
\enq

\nin
\emph{Remark.}
Perhaps surprisingly, even $\xi_\infty$ and the remaining $\xi_s$'s are similarly concentrated.
Since this seems interesting enough 
to record but isn't needed for the rest of what we do (and takes a little while to explain),
we've added its proof as an appendix.

\mn

Let $\xi_\infty^v = \lim_{s \rightarrow \infty} \xi_s^v$. 
Define events  
\beq{events}
\mbox{$Q_v=\{\xi_\infty^v > 0\} \,\, (= \{T_v>r_v\})$ , 
$\,\,R_v = \{r_v \leq M\},\,\,$ and $\,\,Q_v^* = Q_v \cap R_v$.}
\enq
Set $p = 1/L$ ($ = e^{-K}$).

If we cover $W$ in $M$ steps, then $\xi_M \leq (L/\lambda)|\{v \in W : Q_v^* \text{ holds}\}|$; 
so if also $\xi_M \geq m/2$, then 
\beq{eq1}
|\{v \in W : Q_v^* \text{ holds}\}| \geq \gl mp/2.
\enq
So for \eqref{coverW} it is enough to show
\beq{expb}
\pr( (\ref{eq1}) ) = e^{-\Omega(m)}. 
\enq
For the situations we have in mind, this will follow easily from the next two lemmas.
\begin{lemma}\label{assoc}
For any $I\subseteq W$, 
\[
(\pr(\cap_{v\in I} Q_v^*) \leq) \,\,\,\,\,\pr(\cap_{v\in I} Q_v) \leq p^{|I|}.
\]
\end{lemma}
\begin{proof}
Let $\cF_k$ be the $\gs$-field generated by $(X_0,\ldots, X_k)$ and consider the process
\[
S_k^I = \mathbbm{1}_{ \cap_{v\in I} \{T_v > k \wedge r_v\}} \prod_{v\in I} H_v(k \wedge r_v- 1)^{-1}.
\]
We will show that $S_k^I$ is a supermartingale.  
The Martingale Convergence Theorem (Theorem~\ref{conv})
then
says $S_k^I \rightarrow S_\infty^I$ (a.s.) and
\[
\E S_\infty^I \leq \E S_0^I = 1,
\]
which, since 
\[
S_\infty^I = \mathbbm{1}_{\cap_{v\in I} Q_v} \prod_{v\in I} H_v(r_v - 1)^{-1},
\]
gives the desired
\[
\pr(\cap_{v\in I} Q_v) = \E\big[S_\infty^I \prod_{v\in I} H_v(r_v - 1) \big] \leq p^{|I|}.
\]

To see that $S_k^I$ is a supermartingale (here just extending the proof of \cite[Lemma 2]{Aldous}), 
it is enough to show 
\[
\mbox{$\E[S^I_{k+1}|\f_k] \leq S^I_k\,\,$ on $\,\, \{T_v >k\wedge r_v~\forall v\in I\}$}
\]
(since outside this conditioning set, $S_{k+1}=S_k=0$).
But here, with $J=\{v\in I: k< r_v\}$, we have
\begin{align*}
\E(S_{k+1}^I|\cF_k) &= \pr(\cap_{v\in J} \{T_v > k+1\} | \cF_k) \prod_{v\in J} H_v(k)^{-1} 
\prod_{v\in I\sm J} H_v(r_v-1)^{-1} \\
&= (1 - \sum_{v\in J} \pt(X_k, v))\prod_{v\in J} H_v(k)^{-1} 
\prod_{v\in I\sm J} H_v(r_v-1)^{-1} \\
&\leq \prod_{v\in J} (1 - \pt(X_k, v))\prod_{v\in J} H_v(k)^{-1} 
\prod_{v\in I\sm J} H_v(r_v-1)^{-1} \\
&= \prod_{v\in J} H_v(k-1)^{-1} \prod_{v\in I\sm J} H_v(r_v-1)^{-1} = S_k^I.
\end{align*}

\end{proof}

Recalling that $p=e^{-K}$, we now set 
\beq{Kepsgd}
\mbox{$K = \max\{[20e(64+C)]^2, \log (1/\gb)\}$, $\,\, \eps = \gl p/5$,
and $\,\,\, \gd = \eps/K$}
\enq
The reasons for these choices will appear below (see \eqref{middle}), 
and for now we just mention that (i)
the more important constraint in the definition of $K$ is the first, and 
(ii) the main thing to keep in mind here is that there
is nothing preventing us from taking $K$ as large, and $\gd$ as small, as needed
to make things work (cf.\ ``Perspective'' following \eqref{middle});
in particular, the only reason for the fussy specification of $K$ is to make the role of this 
choice a little clearer below.

Define $\vp_\gd(y,z) = \mathbbm{1}_{\{\vp(y,z)\leq \gd\}}\vp(y,z)$ and,
for a multisubset $Y$ of $V$,
\[
\pt_\gd(Y,z) = \sum_{y\in Y }  \pt_\gd(y,z).
\]
(For Theorem~\ref{main} we could skip $\pt_\gd$ and work with $\pt(Y,z)$,
defined in the natural way, but the present version will be needed in Section~\ref{largeB}.)

For the next lemma we take
$\mathcal{Z}$ to be the set of those $Z\subseteq W$ of size at least $2\eps m$ for which
\beq{loc}
\mbox{there is a multisubset $ Y $
of $ V $ of size at most $M$ with 
$\pt_\gd(Y,z) > K/4 \,\,\,\, \forall z \in Z.$}
\enq

\begin{lemma}\label{union} 
There is an $\I\sub \C{W}{\eps m}$ with 
\beq{size}
    |\mathcal{I}| < (15\eps)^{-\eps m}
\enq
such that
\beq{contain}
\mbox{each $Z\in \Z$ contains some $I\in \I$.}
\enq
\end{lemma}

\nin
\emph{Remark.}  
Our eventual bound on the probability in \eqref{expb} will be 
(with an appeal to Lemma~\ref{assoc})
\[
\sum_{I\in \I}\pr(\cap_{v\in I} Q_v^*) \leq |\I|p^{\eps m},
\]
so we want the r.h.s.\ of \eqref{size} to be small relative to $p^{-\eps m}$,
which will be true with the present bound since we will have $\gl\geq 1/2$
(recall $\eps =\gl p/5$).

\begin{proof}[Proof of Lemma \ref{union}]

Fix $Z\in \Z$, let $Y$ be as in \eqref{loc}, and set
\[
W_0 = \{z\in W:\pt_\gd(Y,z) > \sqrt{K} \},
\]
noting that 
\beq{W0}
|W_0| < M/\sqrt{K}.
\enq

Consider the random submultiset $Y'$ of $Y$ gotten by including members of $Y$ 
independently, each with probability $32\gd/K$, and set
\[
N_\gd(Y') = \{z\in W : \pt_\gd(Y',z) \geq \gd\}.   
\]
We assert that with positive probability,  
\beq{Y'size}
|Y'| < 33\gd M/K =:t,
\enq
\beq{Y'hood} 
|N_\gd(Y')\sm W_0|< 64m/\sqrt{K},
\enq
and
\beq{overlap}
|N(Y')\cap Z| > \eps m.
\enq
\emph{Proof.}
Since $|Y'|$ is binomial with parameters $M'\leq M$ and $32\gd/K$, the probability of
violating \eqref{Y'size} is small.

For $z\in W\sm W_0$, we have $\E \pt_\gd(Y',z) \leq 32\gd/\sqrt{K}$, and (by Markov's Inequality)
$\pr(z\in N_\gd(Y'))\leq 32/\sqrt{K}$; so $\E |N_\gd(Y')\sm W_0|\leq32 m/\sqrt{K}$,
and a second application of Markov gives $\pr(\mbox{\eqref{Y'hood} fails})\leq 1/2$.

\nin

Finally, set, for $z\in Z$ and $y\in Y$, $\psi_z=\pt_\gd(Y',z)$ and
$\gz_y=\mathbbm{1}_{\{y\in Y'\}}$.  
Then $\psi_z = \sum\{\gz_y\pt(y,z):\pt(y,z)\leq \gd\}$,
$\E \psi_z > 8\gd$, and $\textrm{Var}(\psi_z) < \gd \E\psi_z$, implying
(e.g.\ by the second moment method; this is reason for the 32)
$\vt:= \max_{z\in Z}\pr (\psi_z< \gd)< 1/6$. 
On the other hand, Markov
gives $\pr(|Z\sm N_\gd(Y')|>3\vt |Z|) < 1/3$,
so $|N_\gd(Y')\cap Z| > (1-3\vt )|Z| > \eps m$ with probability at least 2/3, and
the assertion follows.

\qed

It follows that there \emph{is} some multiset $Y'$ satisfying 
(\ref{Y'size}), (\ref{overlap}), and (from \eqref{W0} and \eqref{Y'hood}; recall $M\leq Cm$)
\[
|N_\gd(Y')| <(64 +C)m/\sqrt{K}.
\]
Thus with
\beq{calY}
\Y = \{Y' \mbox{a multisubset of $V$} : 
|Y'| \leq t,\: |N_\gd(Y')| \leq  (64 +C)m/\sqrt{K} \},
\enq
we find that
\[
\I := \bigcup_{Y'\in \Y}\C{N_\gd(Y')}{\eps m} 
\]
satisfies \eqref{contain}.  But it also satisfies \eqref{size}:

Noting that $|\Y|\leq \C{n+t}{t-1} < \C{2m/\gb}{33C\gd m/K}$ (say), 
and recalling \eqref{Kepsgd} and the bound on $|N_\gd(Y')|$ in \eqref{calY}, we have
\begin{align}
|\cI| &< \C{2m/\gb}{33C\gd m/K}\C{(64 +C)m/\sqrt{K}}{\eps m} \label{|I|}\\
&< \left(\frac{2eK}{33C\gb\gd}\right)^{(33CK^{-2})\eps m}
\left(\frac{e(64+C)}{\eps\sqrt{K}}\right)^{\eps m}\label{middle}
\,\, < \,\,
(15\eps)^{-\eps m}
\end{align}
where we used $\gb> e^{-K}$
to bound the first term in \eqref{middle} by (say)
$e^{(70C/K)\eps m}$.
\end{proof}

\nin
\emph{Perspective.}
There is less here than meets the eye:  the main point is the $1/\sqrt{K}$ of \eqref{Y'hood}, 
which, since \emph{we} choose $K$, can be used to make the second factor on the r.h.s.\ of
\eqref{|I|} much smaller than $\eps^{-\eps m}$;
though we've taken $\gd$ only (roughly) as small as necessary to make the 
first factor irrelevant, there was nothing to stop us from making it smaller, so this factor
was not really an issue; the remaining terms (including the canceling $m$'s) 
may safely be ignored.

\mn
\emph{Proof of Theorem~\ref{main}.}
We prove this with $\gd(C,\gb)$ the $\gd$ of \eqref{Kepsgd},
noting that \eqref{vpvw} then gives $\pt_\gd(v,w)=\pt(v,w)$ 
for relevant $v,w$, 
whence $\gl\approx 1$.
As observed above, we just need \eqref{expb}; namely, with  
$Z = \{w \in W: Q_w^* \text{ holds}\}$,
\beq{prZ}
\pr(|Z| \geq \gl mp/2)< e^{-\gO(m)}.
\enq

If $Q_w^*$ holds (that is, $M\geq T_w> r_w$), then 
$   
\gl/L \geq \prod_{t\leq M} (1-\pt(X_t,w)),
$  
implying (with some room since the $\pt(v,w)$'s are small)
$\sum_{t\leq M} \pt(X_t,w)  > K/2$;
thus $Z$
satisfies (\ref{loc}) (with $Y=\{X_1\dots X_M\}$).
So if the event in \eqref{prZ} holds, then $Z\in \Z$ and we have
$\cap_{v\in I}Q_w^*$ for some $I\in \I$, which according to
Lemma~\ref{assoc} (and \eqref{size}) occurs with probability at most
\[
|\I|p^{\eps m} = e^{-\gO(m)}.
\]\qed

(So here the $K/4$ in \eqref{loc} could have been $K/2$,
but we will need slightly more room in Section~\ref{largeB}.)

\section{Proof of Lemma~\ref{lemmaB}}\label{largeB}

Our main new point here is Claim~\ref{claim1}, given which Lemma~\ref{lemmaB} will 
be another simple application of the material of Section~\ref{mainsec}.
We begin by setting parameters, in particular the $\gD$ of the lemma, noting 
again that 
these fairly careful specifications are meant to make the arithmetic below
easier to track (if one \emph{cares} to track it), 
but that there is nothing delicate in these choices, since there are no constraints on $\gD$
(beyond its being a function of $D$ and $\gc$).  With this advisory, we take
$\gb=\gc/2$, $C=D/\gb$ and $M=Dn $ ($=\gb Cn$); $K$, $\eps$, $\gd$ as in \eqref{Kepsgd}
(again, with $p=e^{-K}$);
\[
\mbox{$d=1/\gd\,\,$ and $\,\, \vr =\gc K p/(160)$;}
\]
and, finally,
\beq{gD}
\gD = 16 Dd^2/(\gc \vr)
\enq
(so $\gD$ is roughly $e^{3K}$).

Set 
\[
\mbox{$B=\{v:d_v> \gD\}\,\,$ and $\,\, S=\{v:d_v\leq d\}$.}
\]

\begin{claim}\label{claim1}
There is a $W\sub B$ of size at least $|B|/2$ such that, with $S^* = N(W)\cap S$,
\beq{eq4.1}
\pr(|\{t\in [M]: X_t\in S^*\}|> \vr n) < e^{-\gO(n)}.
\enq
\end{claim}

\nin
\emph{Proof.}  
We first observe that for all $t$,
\beq{3.1item}
\sum_{v\in B}\pt^t(v,S) \leq (d/\gD)|S|.
\enq

\begin{proof} 
Using reversibility (which implies $\pi (v)\pt^t(v,w)=\pi(w)\pt^t(w,v)$; see \eqref{pathrev}), we have
\[
\sum_{v\in B}\pt^t(v,S)  =\sum_{v\in B}\sum_{w\in S}\pt^t(v,w)
= \sum_{w\in S}\sum_{v\in B}\pt^t(w,v)d_w/d_v
\leq (d/\gD)\sum_{w\in S}\pt^t(w,B) \leq (d/\gD)|S|
\]
\end{proof}

It follows that for all $T$ (now just using $|S|\leq n$),
\[
|B|^{-1} \sum_{v\in B}\pr_v(X_{[T]}\cap S\neq \0) \leq Td/(\gc\gD)
\]
and 
\[
|\{v\in B: \pr_v(X_{[T]}\cap S\neq \0) \geq 2Td/(\gc\gD)\}| \leq |B|/2;
\]
so if we set $T=\gc\gD/(4d)$ 
and take
\[
W=\{v\in B: \pr_v(X_{[T]}\cap S\neq \0) < 2Td/(\gc\gD) \,\, (=1/2)\},
\]
then 
\[
|W|\geq|B|/2
\]
and we just have to show 
\beq{Wsatisfies}
\mbox{$W$ satisfies \eqref{eq4.1}.}
\enq

To see this, let $\xi_i$ be the time between the $(i-1)$st and $i$th visits to $S^*$.
Then (independent of history up to the $(i-1)$st visit),
\beq{prxii}
\pr(\xi_i> T)> 1/(2d).
\enq
(Starting from $v\in S^*$, we're in $W$ at the first step with 
probability at least $1/d$ and then with probability at least 1/2 the time to return
to $S^*$ is at least $T$.) 
With
$\psi_i =  \mathbbm{1}_{\{\xi_i>T\}}$ and
$\psi=\sum_{i=1}^{\vr n} \psi_i$,
visiting $S^*$ more than $\vr n$ times 
(the event in \eqref{eq4.1}) requires 
\[
\psi <M/T = \vr n/(4d).
\]
But $\psi$ stochastically dominates $\psi'\sim \textrm{Bin}(\vr n,1/(2d))$
(by \eqref{prxii}), and 
$\pr(\psi'< \vr n/(4d)) < e^{-\gO(n)}$.

This completes the proofs of \eqref{Wsatisfies} and Claim~\ref{claim1}.
\qed

\mn
\emph{Proof of Lemma~\ref{lemmaB}.}
We use the machinery of Section~\ref{mainsec} with $W$ as in Claim~\ref{claim1}
(and $|W|=m$) and other parameters as in the first paragraph of this section.

[One picky adjustment:  
If $v$ is pendant (i.e.\ of degree one) with unique neighbor $w$,
then $\pt(v,w)=1$ and the $\gl$ of \eqref{lambda} can be zero.
But, except when $X_0=v$, transitions from pendant vertices 
have no effect on anything in Section~\ref{mainsec}, since
the mandatory next vertex has already been seen and is
no longer contributing to the martingale.  So for the present application we may
without penalty modify \eqref{lambda} to require $d_G(v)\geq 2$
(and---getting sillier---exclude the unique neighbor of $X_0$ from $W$ if $X_0$ happens
to be pendant); thus
we assume for this little discussion that $\gl\geq 1/2$.]

Define events
\[
E= \{|\{v\in W :\mbox{$Q_v^*$ holds}\}| > \gl p m/2\}
\]
and
\[
F= \{|\{t\in [M] : X_t\in S^*\}| \leq \vr n\}
\]
($S^*$ as in Claim~\ref{claim1}).
As earlier, we just need to show \eqref{expb} (namely, $\pr(E) < e^{-\gO(n)}$), 
which in view of \eqref{eq4.1} will follow from
\[
\pr(E\wedge F) < e^{-\gO(n)}.
\]

To see this, note first that, with
\[
Z_0 = \left\{v:\sum\{\pt(X_t,v):t\in [M],  X_t\in S^*\}> K/4\right\},
\]
$F$ implies 
$|Z_0|< 4\vr n/K$. 

As in the proof of Theorem~\ref{main}, if $Q_v^*$ holds, then
$ 
\gl/L \geq \prod_{t\leq M} (1-\pt(X_t,v)),
$ 
which in the present situation (i.e.\ where $\gl \geq 1/2$) implies 
\[
\sum_{t\leq M} \pt(X_t,v)  \geq (2\log 2)^{-1}\log (L/\gl)  > K/2;
\]
so if also $v\in W\sm Z_0$, then (since $d=1/\gd$),
\[
\sum_{t\leq M} \pt_\gd(X_t,v)  > K/4.
\]

So if $E\wedge F$ holds, then
\[
Z:= \{v\in W:\mbox{$Q_v^*$ holds}\}\sm Z_0
\]
satisfies \eqref{loc} with
\[
Y=\{X_t: t\in [M], X_t\not\in S^*\},
\]
and (with minor arithmetic, again using $\gl\geq 1/2$)
\[
|Z| >\gl p m/2 - 4\vr n/K > 2\eps m;
\]
that is, $Z\in \Z$.  We thus have
$\cap_{v\in I}Q_v^*$ for some $I\in \I$, and, by
Lemma \ref{assoc}, 
\[
\pr(E\wedge F) \leq |\I|p^{\eps m} = e^{-\gO(m)}.
\]\qed

\section{Partitions}\label{Partitions}

Here we prove Theorem~\ref{thmMT}.  As mentioned earlier, this will be based on
Corollary~\ref{corP} provided we exclude two possibilities---\eqref{LBineq} and 
\eqref{LYineq}---that imply the conclusion of the theorem for other reasons
(as shown earlier in Lemma~\ref{lemmaB} and soon in Lemma~\ref{lemmaY}).

We fix $C$ and consider a walk of length $Cn$ on the ($n$-vertex) graph $G$.
Let $\gd =(2/3)\gd(C/\gc^2,\gc)$
(see Theorem~\ref{main} for $\gd(\cdot,\cdot)$, Corollary~\ref{corP}
for our intended use, and \eqref{3gd2}
for the silly reason for the 2/3), 
and let $\gD$ be as in Lemma~\ref{lemmaB} with $C$ in place of $D$, and
\[    
\theta = \gd^2.
\]   
(This extra parameter could be skipped,
but is included as it will appear pretty often.)

Set (for any $v$ and $R$)  
\[
\mbox{$B_v (R)= \{w\neq v:\pr_v(w\in X_{[R]})> \gd/2\}$}
\]
and
\[
\mbox{$B'_v(R) = \{w\neq v:\pr_v(T_w\leq \min\{R,T_v^+\})> \gd/2\}$}.
\]

\mn
(We don't actually need the superset $B_v(R)$ of $B_v'(R)$,
but keep it to point out that the upper bound shown in Lemma~\ref{lemmaQ}
doesn't use the extra constraint in $B'_v(R)$.)

\mn
\emph{Preview.}  
For any specification of $V_i$'s we will take 
\beq{Wi}
U_i = \{w\in V_i:\max_{w\neq v\in V_i}\pt_i(v,w)<\gd\},
\enq
Thinking of $v$'s that cause exclusions from these $U_i$'s, we 
say $v\in W\sub V$ is \emph{good for W} 
(or just \emph{good} if the identity of $W$ is clear) 
if
\beq{good}
\max_{v\neq w\in W}\pt_W(v,w)<\gd.
\enq
We are hoping for $V_i$'s in which few vertices are bad (not good), in which case 
we can use the trivial
\beq{Uimost}
|V_i\sm U_i|\leq |\{v\in V_i:\mbox{$v$ is bad for $V_i$}\}|/\gd.
\enq

Perhaps surprisingly---and luckily, since other options seem difficult---much 
of our production of such $V_i$'s (all but what's covered by Lemma~\ref{lemmaR})
can be based on the following easy point.

\begin{obs}\label{Bvobs}
For any R, sufficient conditions for $v$ to be good for W are
\beq{good1}
W\cap B_v'(R)=\0
\enq
and
\beq{good2} 
\pr_v(X_{[R]}\cap W=\0)< \gd/2.
\enq
(These are enough since then for $w\in W\sm \{v\}$, 
\[
\pt_W(v,w) \leq \pr_v (T_w\leq \min\{R,T_v^+\}) + \pr_v(X_{[R]}\cap W =\0) < \gd.)
\]
\end{obs}

\nin
\emph{Note.}  As mentioned earlier, a reader interested in a warm-up for what we're about to do
might find this a good time to take a look at Section~\ref{expand}.


\mn

Before turning to our main line of argument we 
dispose of an easy case. 
Say $v$ is $(\gd,R)$-\emph{recurrent} if 
\beq{vinR'}
\pr_v(T_v^+\leq R) >1-\gd
\enq 
and \emph{$(\gd,R)$-transient} otherwise.

\mn
\begin{lemma}\label{lemmaR}
If, for some $R$,
\beq{recurrents}
|\{v: \mbox{$v$ is $(\gd,R)$-recurrent}\}| > 2\gc n,
\enq
then $G$ admits a partition as in Corollary~\ref{corP} with $\vt = \gd\gc/(2R)$.
\end{lemma}

\mn
\emph{Proof.}
Let $S$ be the set in \eqref{recurrents}.  Notice that, for any $v,w$,
\beq{BvBw}
\mbox{if $d_G(w)\geq d_G(v)$ and $v\in B_w'(R)$, then $w\in B_v'(R)$}
\enq
(since 
$\pr_v(T_w\leq \min\{R,T_v^+\}) = (d_G(w)/d_G(v))\pr_w(T_v\leq \min\{R,T_w^+\}) $;
see \eqref{pathrev}).

Let $\gG$ be the graph on $S$ with $v\sim w$ if $w\in B'_v(R)$ or vice versa.
Order $V$ by some ``$\prec$'' with $v\prec w \Ra d_G(v)\leq d_G(w)$ and notice that
\eqref{BvBw} implies (the first inequality in)
\[
d^+(v) \leq |B_v'(R)| < 2R/\gd  \,\,\,\, \forall v\in S
\]
(where $d^+(v)=|\{w: v\prec w\sim v\}|$), whence the chromatic number of
$G$ is at most $2R/\gd$ (a standard exercise or see e.g.\ \cite[Section 5.1]{MGT}).
We now take $\{W_j\}$ to be a (proper) $(2 R/\gd)$-coloring of $\gG$, 
and notice that, for any $j$ and distinct $v,w\in W_j$,
\beq{3gd2}
\pt_{_{W_j}}(v,w)\leq \pr_v(T_w\leq \min\{R,T_v^+\}) + \pr_v(T_v^+> R) < 3\gd/2.
\enq
(as in \eqref{Ui}).  We also have (with $\vt$ as in Lemma~\ref{lemmaR}) 
\[
\sum\{|W_j|:|W_j|\leq \vt n\}\leq (2R/\gd)\vt n =\gc n;
\]
so we satisfy the demands of Corollary~\ref{corP} by taking
$\{V_i\}= \{W_j:|W_j|> \vt n\}$ and $U_i=V_i$ $\forall i$
(and $V^0=V\sm \cup V_i$).

(For clarity we just note that the bound we actually need in \eqref{3gd2} is 
$\gd(C/\gc,1)> 3\gd /2 $.)\qed

\mn

We now turn to the main argument.
Fix $k$ with
\[
\gd^{k-3} < (16\sqrt{\gD})^{-1}
\]
($k=5$ will do since $\gD$ is roughly $\gd^{-3}$), and
let $N$ be minimum with
\beq{1-gdM}
(1-\gd)^N < \gd^k.
\enq

Let $R_0=1$ and, for $i\geq 1$,
\[
R_i = 4CN\gd^{-10} R_{i-1}.
\]
Choose $i< 10\gd^{-k}$ for which at least $.9n$ vertices $v$ satisfy
\beq{RR'}
\pr_v(T_v^+\in (R_{i-1},R_i])<\gd^k.
\enq
Parameters we will use (collected here to have them in one place,
though it will take us a little while to get to the $Q_i$'s) are then:
\[
\mbox{$R'=NR_{i-1}$, $\,\,Q=(4/\gd)R'$, $\,\,R=R_i$,}
\]
\[
\mbox{$Q_1 = Q\gd^{-3}, \,\,$
  and \,
$Q_2 = Q_1\gd^{-2}\theta^{-1}=Q_1\gd^{-4}$}
\]
(so $R = CQ_2/\theta$).
We also abbreviate 
\[(R',R] = I,\]
since this interval will appear frequently.
(The ratios between parameters are generous but convenient,
in particular supporting occasional use of inequalities of the form 
$e^{-1/\gd}< \gd^{O(1)}$, which hold since $\gd$ is small.)
For minor reasons at \eqref{pick} we want---and, to avoid very silly distractions, 
will just assume---
\beq{R'+1}
\mbox{$R'+1$ is even.}
\enq

In view of Lemmas~\ref{lemmaB} and \ref{lemmaR}, we may assume
\beq{T}
\mbox{at least $.6n$ vertices $v$ are $(\gd, R)$-transient, have $d_G(v) < \gD$, 
and satisfy \eqref{RR'}.}
\enq
Let $\T$ be the set of such $v$'s.

\mn
\begin{lemma}\label{lemmaQ}  For any $v\in \T$, 
$\,\,\,
|B_v(R)\cap \T|\leq Q.
$
\end{lemma}

\nin
\emph{Proof.}
We first observe that 
\[
\pr_v(v\in X_I) 
< 2\gd^{k-1}.
\]

\nin
[Because:  If $\{t\geq 0:X_t=v\} =\{t_0<t_1< t_2<\cdots\}$, 
then
$v\in X_I$ implies that either 
\beq{either}
t_i-t_{i-1}\leq R_{i-1} \,\,\,\,\forall i\in [N]
\enq
or, for some $j\leq N$,
\beq{or}
\mbox{$t_i-t_{i-1}\leq R_{i-1}~$ for $~ i\in [j-1]~$ and $~t_{j-1}+R_{i-1}<t_j \leq R_i$.}
\enq
But by \eqref{1-gdM} and \eqref{RR'}
the probabilities of \eqref{either} and \eqref{or} are
less than (respectively) $\gd^k$ and
\[
\gd^k\sum_{j\in [N]}(1-\gd)^{j-1} < \gd^{k-1}.]
\]

Set $\ell_v(I)= |\{t\in I: X_t=v\}|$ 
and notice that, for any $v\in \T$,
\beq{Evlv}
\E_v \ell_v(I) < 2\delta^{k-2}
\enq
(since 
$
\E_v\ell_v(I) =\sum_{u\geq 1}\pr_v(\ell_v(I)\geq u)\leq 2\gd^{k-1}\sum_{u\geq 1} (1-\gd)^{u-1} = 2\gd^{k-2}).
$ 

\mn
It follows that for distinct $v,w \in \T$,
\begin{align}
(\pr_v(w\in X_I)\leq)\,\,\,\,\,     \E_v \ell_w(I) &= \sum_{t \in I} \pt^t(v,w) \nonumber\\
    &\leq\label{AF} \sqrt{\gD}\left[\sum_{t\in I} \pt^{2\lfloor t/2 \rfloor}(v,v)\sum_{t\in I}\pt^{2\lceil t/2 \rceil}(w,w)\right]^{1/2}\\
    &\leq\label{pick}  \sqrt{\gD}\left[2\sum_{t \in I} \pt^t(v,v)2\sum_{t \in I} \pt^t(w,w)\right]^{1/2}\\
    &<\label{bound} 4 \sqrt{\gD}\cdot \delta^{k-2}  < \delta/4,
\end{align}
where \eqref{AF} is Lemma~\ref{AF320} and Cauchy-Schwarz, 
\eqref{pick} uses Lemma \ref{mono} and \eqref{R'+1},
and \eqref{bound} is given by \eqref{Evlv}.

Thus, finally,
\[
B_v(R)\cap \T \sub \{w \in \T\sm \{v\} : \pr_v(w \in X_{[R']}) > \delta/4\},
\]
a set of size at most $(4/\gd)R' =Q$.

\qed

\nin
\emph{More preview.}
In what follows, aiming for Corollary~\ref{corP}, we will discard $V\sm \T$ 
(that is, include it in $V^0$) and consider a random 
partition of $\T$, hoping to use
Observation~\ref{Bvobs} (and the discussion preceding it) to say that (with good probability) 
much of $\T$ lies in blocks that behave as the corollary requires.
Roughly speaking, what we get from Lemma~\ref{lemmaQ} is likelihood 
of \eqref{good1}:  if the number of blocks in our random partition
is much larger than $Q$, then the block containing $v$ is unlikely to meet $B_v'(R)$.

For \eqref{good2}
a natural intuition is that ``transience'' (failure of \eqref{vinR'}) implies 
that, for the walk started from $v$, $X_{[R]}$ is likely to be large, which,
suitably quantified, does imply that \eqref{good2} is likely (for $v$ and its random block $W$).
This intuition turns out to be not quite correct, but, as shown in Lemma~\ref{lemmaY},
if it is wrong too often then the conclusion of Theorem~\ref{thmMT}
holds for other (simpler) reasons.

Set 
\[
\D =\{v\in \T: \pr_v(|X_{[R]}\cap \T|<Q_1) > \theta\gd\}.
\]

\begin{lemma}\label{lemmaY}
If
\beq{LYineq}
|\D|\geq 2\theta n,
\enq
then $\pr(X_{[Cn]}\supseteq V) =e^{-\gO (n)}$.
\end{lemma}

\nin
\emph{Proof.}
We first claim that
\beq{anyv}
\mbox{for any $v\in V$,  $\,\,\,\pr_v(|X_{[R]} \cap \D| > Q_2) < \exp [-1/\gd].$}
\enq

\nin
\emph{Proof.}
With $\{X_t\}$ started from $v$, let $t_0=\min\{t:X_t\in \D\}$ and, for $i\geq 1$,
\[
t_i=\min\{t:  X_t\in \D, \,\, |X_{(t_{i-1},t]}\cap\T|\geq Q_1\}.
\]
(That is, $t_i$ is the first time that the walk is in $\D$,
having seen at least $Q_1$ distinct vertices of $\T$ since $t_{i-1}$.)

For the event in \eqref{anyv} we must have (\emph{very} generously) 
\[
\mbox{$t_i-t_{i-1}\leq R\,\,$   $\forall i\in [Q_2/Q_1]$,}
\]
which, since each $X_{t_{i-1}}$ is in $\D$, occurs with probably less than 
$  
(1-\theta \gd)^{Q_2/Q_1} < e^{-1/\gd}.
$   

\qed

We can now show
\beq{tinCn}
\pr(|X_{[Cn]}\cap \D| \geq 2 CQ_2 n/R) < e^{-\gO(n)},
\enq
which gives the lemma since $CQ_2/R=\theta$.

\mn
\emph{Proof of \eqref{tinCn}.}
For $i\in [Cn/R]$ let $\xi_i$ be the indicator of
\[
\{|X_{((i-1)R,iR]}\cap \D|>Q_2\}.
\]
Then
\[
 |X_{[Cn]}\cap \D| 
\leq R\sum \xi_i + CQ_2 n/R,
\]

\mn
so the event in \eqref{tinCn} requires $\xi:=\sum\xi_i> CQ_2n/R^2$.
But $\xi$ is stochastically dominated by 
$\xi'\sim \textrm{Bin}(Cn/R,e^{-1/\gd})$
(by\eqref{anyv}), and 
$\pr(\xi'> CQ_2n/R^2) < e^{-\gO(n)}$.

\qed

So we may assume 
\beq{Dsmall}
|\D|<2\theta  n.
\enq

For the partition of Corollary~\ref{corP}, we include
$V\sm \T$ in $V^0$ and will mainly be interested in $\T\sm \D$.  
Setting 
\[
\gz = \theta/Q,
\] 
we randomly (uniformly)
partition $\T$ into $\gz^{-1}$ blocks, 
usually called $W$, 
and want to say that each $v\in \T\sm \D$ is likely to be good 
(meaning, of course, good in its block).

\mn
\begin{lemma}\label{lemmaZ}
If $v\in \T\sm\D$ then $\pr(\mbox{$v$ bad})  <4\theta$.
\end{lemma}

\nin
\emph{Proof.}  
We want to say that, at least for $v\in \T\sm \D$, 
\eqref{good1} and \eqref{good2} are likely for $v$ and 
the block $W$ containing it.
For \eqref{good1} this is just 
\beq{PwBv}
\pr(W\cap B_v'(R)\neq \0) < \gz |B_v'(R)\cap\T| < \gz Q = \theta
\enq
(this just requires $v\in \T$;  see Lemma~\ref{lemmaQ}).

For \eqref{good2} (now using $v\not\in \D$),
with unsubscripted $\pr$ referring to the choice of the block $W$ containing
$v$ \emph{and} the walk from $v$, we have
\begin{eqnarray}
\E_W[\pr_v(X_{[R]}\cap W=\0) ]
&=&   \pr (X_{[R]}\cap W=\0)  \label{prwV}\\
&< & \pr_v(|X_{[R]}\cap \T| < Q_1) + e^{-\gz Q_1} 
< \theta\gd+ e^{-1/\gd}=:q .\nonumber
\end{eqnarray}
But the l.h.s.\ of \eqref{prwV} is at least
\[
(\gd/2)\pr_W[\pr_v(X_{[R]}\cap W=\0)\geq \gd/2 ],
\]
so 
\beq{PwPv}
\pr_W[\pr_v(X_{[R]}\cap W=\0)\geq \gd/2 ] < 2q/\gd < 3\theta.
\enq

Combining \eqref{PwPv}
and \eqref{PwBv} now completes the proof of Lemma~\ref{lemmaZ}.

\qed

Again considering our random partition, and using \eqref{Dsmall}
and Lemma~\ref{lemmaZ}, we find that there \emph{exists} a partition
$\{W_i:i\in [\gz^{-1}]\}$ of $\T$ with (say)
\beq{|W|}
|W_i| > \gz n/2 \,\,\forall i
\enq
and 
\[
|\{v:\mbox{$v$ bad}\}| < 5\theta|\T\sm \D| +|\D| < 7\theta n
\]
(where, again, ``$v$ bad'' means bad in its $W_i$).

Say $W_i$ is \emph{nice} if 
\[
|\{v\in W_i: \mbox{$v$ bad}\}| < \gd |W_i|/2,
\]
noting that this implies 
\[
|U_i| > |W_i|/2
\]
(recalling that $U_i$ was defined in \eqref{Wi} and using \eqref{Uimost}).

On the other hand,
\[
\sum\{|W_i|:\mbox{$W_i$ not nice}\} \leq (2/\gd)|\{v:\mbox{$v$ bad}\}|
< 14\theta n/\gd, 
\]
whence $\sum\{|W_i|:\mbox{$W_i$ nice}\} >|\T|-14\theta n/\gd > .5 n$
(see \eqref{T});
so, with $\vt =\gz/2$ (see \eqref{|W|}), the collection $\{V_j\}$ of nice $W_i$'s,
with $V^0=V\sm \cup V_j$, is the desired partition.

\qed

\section{Expanders}\label{expand}

As promised near the end of Section~\ref{Intro}, 
this and the next section give separate treatment to
two previously known cases of 
Theorem~\ref{thmMT}, as relatively simple illustrations of the use of 
Corollary~\ref{corP}.
Here we provide (a little sketchily) a simpler substitute for much of
Section~\ref{Partitions} in the case of expanders (for which, as said earlier,
Theorem~\ref{thmMT} was observed in \cite{Yehudayoff} to follow easily 
from \cite{Gillman}).
Note we are still using the defaults $G=(V,E)$ and $|V|=n$.

Suppose the transition matrix, $P$, of RW on $G$ has eigenvalues 
$1 = \gl_1 \geq \cdots \geq \gl_n \geq -1$ (as guaranteed by 
Perron-Frobenius).
We call $G$ an \textit{$\eps$-expander} if
$\max\{|\gl_2|, |\gl_n|\} < 1-\eps$.
We should show:
\begin{thm}\label{exp}
For RW $(X_t)$ on an $\eps$-expander $G$,
\[
\pr(X_{[Cn]} =V) = \exp[-\gO_{\eps,C}|V|].
\]
\end{thm}
(Note $\gd$, $\gD$ are still as in the second paragraph of Section~\ref{Partitions}.)
In view of Lemma~\ref{lemmaB}, we may assume at least $(1-\gc)n$ 
vertices of $G$ have degree at most $\gD$. Let $\T$ be the set of such vertices. 
Application of Observation~\ref{Bvobs} here will be based on the next two assertions.

\begin{prop}\cite[Theorem 5.1]{Lovasz}\label{tail}
For an $\eps$-expander G and $S\subseteq V$,
\[
|\pt^t(v,S) - \pi_S| \leq \sqrt{\pi_S/\pi_v}\,(1-\eps)^t.
\]
\end{prop}
\begin{prop}\label{hitting}
For RW
on an $\eps$-expander $G$ and $S\subseteq V$,
$\,\,
\pr(T_S > t)  < \left(1-\pi_S/2\right)^{\eps t/(2\log n)}.
$
\end{prop}
\nin
[We include the trivial proof:  
Set $s=2\log n/\eps$.  Proposition~\ref{tail} gives (say) 
$\pr(X_{r+s} \in S|X_r=v) > \pi_S/2$ for any $r$ and $v$, so
\[
\pr(T_S >t) \leq \pr(X_{ks}\not\in S ~\forall k\in [t/s]) < \left(1-\pi_S/2\right)^{\eps t/(2\log n)}. ]
\]

Now thinking of \eqref{good1}, we observe that 
there is a fixed $Q$ such that for any $R=o(n)$ and $v$,
\beq{BBQ}
(|B_v'(R)\cap \T|\leq ) \,\,\,\,|B_v(R)\cap \T| <Q.
\enq
[Because:  By Proposition~\ref{tail}, there is a fixed $T$ (depending on $\eps,\gd,\gD$)
so that, for any $w\in \T$,
\[
\pr_v(w\in X_{(T,R]}) < \sqrt{\gD}\eps^{-1}(1-\eps)^T + R\pi_w < \gd/4;
\]
so $B_v(R)\cap \T\sub \{w:\pr_v(w\in X_{[T]})>\gd/4\}$, a set of size less than $4T/\gd=:Q$.]

On the other hand, Proposition~\ref{hitting} \emph{guarantees} \eqref{good2}
whenever $R=\go(\log n)$ and $|W|=\gO(n)$.

Now set $R=\sqrt{n}$ (we need $\log n \ll R\ll n$) and $\gz=\theta/Q$
(recall $\theta =\gd^2$), and consider a random (uniform) partition, $\{W_i\}$, of $\T$ into
$\gz^{-1}$ blocks.  By Observation~\ref{Bvobs} and the discussion above, the probability
that $v\in \T$ is bad in its block $W$ is less than
\[
\pr(W\cap B_v(R))=\0) + \pr(|W|< \gz n/2) < \gz Q +o(1) =\theta +o(1).
\]

The rest of this is essentially the same as the end of Section~\ref{Partitions}
(following the proof of Lemma~\ref{lemmaZ} and omitting $\D$); so we 
won't duplicate, but briefly:
The preceding discussion shows existence of a partition $\{W_i\}$ of $\T$ with (say)
$|W_i| > \gz n/2=:\vt n$ $\forall i$, and only $2\theta n$ bad $v$'s.
We then discard (add to $V\sm \T$ to form $V^0$) any $W_i$'s that are 
``not nice,'' meaning 
$|\{v\in W_i: \mbox{$v$ bad}\}| > \gd |W_i|/2$, 
and take $\{V_j\} = \{\mbox{nice $W_i$'s}\}$.

(The definition of ``nice'' is chosen so that $W_i$ nice implies $|U_i|> |W_i|/2$
(see \eqref{Wi} for $U_i$), and the bound on the number of bad $v$'s,
with $\theta\ll \gd$, implies that the number of discarded vertices is small.)

\qed

\nin
\emph{Remark.}  This could also have been handled deterministically, as in the proof of 
Lemma~\ref{lemmaR}, but the intention here was to parallel the main argument of
Section~\ref{Partitions}.

\section{Trees}\label{Trees}

Here we give the promised alternate proof of Theorem~\ref{thmMT} for trees. 
This is again based on Corollary~\ref{corP}, but now without Observation~\ref{Bvobs}. 
The proof is constructive (unlike that of Section~\ref{Partitions}) and gives more
than the corollary requires:
\begin{thm}\label{TreePart}
For RW on a tree T, $\gd>0$ and $t=1/\gd$, 
there is a partition $V=W_1\cup\cdots\cup W_k$
with $k \leq (t+1)t^{t+1}$ and (for all $i$)
\[
\max\{\pt_{W_i}(v,w):  v,w\in W_i, v\neq w\} \leq \gd.
\]
\end{thm}

\nin
(To get a partition as in Corollary~\ref{corP} from this, set
$\vt = (2k)^{-1}$,
and take $\{V_j\} = \{W_i : |W_i|\geq \vt n\}$, $U_i = V_i$, 
and $V^0 = V\sm \cup V_i$, noting that $|V^0|\leq n/2$.)

\mn

Our construction is based especially on the following easy property of trees (see e.g.\
\cite[Prop.\ 2.3]{Lovasz}), in which $d(\cdot,\cdot)$ is distance.
\begin{prop}\label{far}
For distinct vertices v, w of T,
$\,\,\pr_v(T_w < T_v^+) \leq 1/d(v,w).$
\end{prop}

\nin
\emph{Usage.}
We regard trees as rooted.
As usual, $v$ is an \emph{ancestor} of $w$ (and $w$ a \emph{descendant} of $v$) if
$v$ lies on the path joining $w$ to the root.
We use $D_v$ for the set of descendants of $v$, $v\wedge w$ for the 
most recent common ancestor of $v$, $w$ (the one furthest from the root), 
and $L_i$ for the set of vertices at distance $i$ from the root.

We will find it convenient to treat partitions as colorings (of $V$).  
We say $W\sub V$ is \emph{safe} if 
\[
\max\{\pt_W(v,w):  v,w\in W, v\neq w\} \leq \gd,
\] 
and a coloring $\gs$ is safe if $\gs^{-1}(c)$ is safe for every $c$.
Since (trivially) $\vp_W(v,w)\leq \vp_{W'}(v,w)$ whenever $v,w\in W\sub W'$, 
Proposition~\ref{far} implies
\beq{redstar} 
\mbox{if $W_1, \ldots$ are safe and $d(W_i,W_j)\geq 1/\gd$ $\,\forall i\neq j$,
then $\cup W_i$ is safe.}
\enq

For the partition of Theorem~\ref{TreePart}
the main thing we have to show is:
\begin{claim}\label{Lt}
For any $T$, there is a safe coloring of $L_t$ with at most $(t+1)t^t$ colors.
\end{claim}
\begin{proof}[Proof of Theorem~\ref{TreePart} given Claim~\ref{Lt}]
Let $\D_q$, $q\in [t]$, be disjoint sets of colors, each of size $(t+1)t^t$.
By \eqref{redstar} it is enough to find, for each $q$ and $i\equiv q\pmod{t}$,
a safe coloring of $L_i$ using colors from $\D_q$.
For $i\geq t$ this is accomplished by applying Claim~\ref{Lt} to $D_v\cap L_i$ for each
$v\in L_{i-t}$ (and again using \eqref{redstar});
for smaller $i$, we can apply the claim to the tree gotten from $T$ by adding a new root
and a path of length $t-i$ joining it to the root of $T$.
(Or check that the proof of the claim also applies here.)
\end{proof}

\begin{proof}[Proof of Claim~\ref{Lt}]
Let $\cee_B,\cee_1\dots \cee_t$ be disjoint sets of colors of size $t^t$.
We color $L_t$ in stages.  For a given stage, we use $U$ for the set of uncolored vertices
at the beginning of the stage, and, for 
$v\in L_0\cup \ldots \cup L_{t-1}$, $U_v = D_v \cap U$.
The process continues until $|U|\leq t^t$, at which point we complete the coloring by
assigning distinct colors from $\cee_B$ to the vertices of $U$.

If $|U|>t^t$, we choose $v\in L_i$ with $|U_v|> t^{t-i}$ and $i$ as large as possible
(so $|U_w|\leq t^{t-j}$ for each $j$ and $w\in D_v \cap L_j$). 
Call $S\sub U_v$ \emph{primitive} (w.r.t.\ $v$)
if $w\wedge z =v$ for all distinct $w,z\in S$.
For $j=1,\ldots ,$ let $S_j$ 
be a maximal primitive subset of 
$U_v \sm (S_1 \cup \ldots \cup S_{j-1})$, ending, say at $S_\ell$,
as soon as the largest surviving primitive set has size less than $t$. 
Thus each of $|S_1|\dots |S_\ell|$ is at least $t$ and, by our choice of $i$,
\[
\ell \,\, (\leq \max\{|U_w|:\mbox{$w$ a child of $v$}\}) \,\, \leq t^{t-i-1};
\]
so we may assign $S_1\dots S_\ell$ distinct colors from $\cee_i$
(and could have taken $|\cee_i|=t^{t-i-1}$).
This completes the stage and leaves $v$ with fewer than $t^{t-i}$ uncolored descendants
(since fewer than $t$ of its children now have $\emph{any}$ uncolored descendants.
Since each $v$ is ``processed'' at most once, we eventually have $|U|\leq t^t$ and
(as above) finish the coloring using $\cee_B$.  

It remains to show that the coloring, $\gs$, is safe.  
Suppose instead that $\gs_w=\gs_z =c$ (for some $w\neq z $ and $c$).
Since $|\gs^{-1}(c)|\leq 1$ for $c\in \cee_B$, we have $c\in \cee_i$
for some $i$.
But then (e.g.) $w$ was colored as part of a primitive set 
$S=\{w_1\dots w_s\}$, with $s\geq t$ and common ancestor $v\in L_i$;
so, since the path from $z$ to $w$ includes $v$, we have
$\vp_c(z,w_j)\geq \vp_c(z,w)$ $\forall j$
(with equality if $z\neq w_j$), where $\vp_c=\vp_W$ with $W=\gs^{-1}(c)$.
Thus
$\vp_c(z,w)\leq 1/s\leq \gd$.
\end{proof}

\nin
\textbf{Acknowledgments.}  We thank Bhargav Narayanan for helpful conversations
and Ori Gurel-Gurevich for telling us the problem, long ago.

\section{Appendix: concentration}\label{appendix}

Usage here is as in Section~\ref{mainsec}, and $v$ will always be a vertex of $ W$.  
As promised following \eqref{H-A},
we show that each $\xi_s$ ($=\xi^W_s$)
is exponentially concentrated about its mean.
\begin{thm}\label{alltime}
For any $\vt >0$,
\[\pr(|\xi_s - m| > \vt m) \leq 2e^{-\vt^2 \gl^{4} m/(8L)^2}.\]
\end{thm}
\nin 
Since $Q_v = \{\xi_\infty^v > 0\}$, 
this gives exponential tail bounds for $|\{v : Q_v\}|$. 
(This isn't quite concentration about the mean since we only know 
$|\{v : Q_v\}|L \leq \xi_\infty \leq (L/\gl)|\{v : Q_v\}|$.)

Theorem~\ref{alltime} is proved using a better
martingale analysis, based on an idea from \cite{aglc}. 
We set $Z_i = \xi_i-\xi_{i-1}$ and  $Z = \sum_{i=1}^s Z_i$ ($=\xi_s-\xi_0$),
and as usual want to bound 
$\E[e^{\gz Z}] $ (with $\gz>0$ to be specified).
The main point here, an instance of \cite[Lemma 3.4]{aglc}, is that we can replace the usual
product of worst case bounds in 
\[
\E[e^{\gz Z}] \leq 
\prod_{i=1}^s \max_{H_{i-1}}\E[e^{\gz Z_i } \mid H_{i-1}] 
\]
by a worst case product:
\begin{lemma}\label{adaptive}
With each $H_i$ ranging  over events $\{X_0 = x_0, X_1 = x_1,\ldots, X_i = x_i\}$, 
\beq{eq.adaptive}
\E[e^{\gz Z}] \leq 
\max\{\prod_{i=1}^s \E[e^{\gz Z_i } \mid H_{i-1}] : H_0 \supseteq H_1 \supseteq \cdots \supseteq H_{s-1}\}.
\enq
\end{lemma}
The next observation will be used to bound the factors in \eqref{eq.adaptive}.
\begin{prop}\label{expbound}\cite[Proposition 3.8]{aglc}
Suppose the $\Re$-valued random variable $Y$ with $\E[Y] = 0$ satisfies 
\[|Y| \leq c\]
and
\[\E[|Y|] \leq M.\]
Then for $|\gz| c \leq 1$,
\[\E[e^{\gz Y}] \leq e^{8 \gz^2Mc}.\]
\end{prop}

Let $H_i = \{X_0 = x_0, X_1 = x_1,\ldots, X_i = x_i\}$ 
(as in Lemma~\ref{adaptive}).
Applying Proposition~\ref{expbound} to each $Z_i \mid H_{i-1}$, with $c= L/\gl^{2}$ 
(see \eqref{xisxis})
and $M=M_i := \E[|Z_i\mid H_{i-1}|]$, gives
\beq{product}
\prod_{i=1}^s \E[e^{\gz Z_i} \mid H_{i-1}] \leq e^{8\gz^2 (L/\gl^2)\sum M_i}
\,\,\,\,\, \mbox{for $|\gz|c \leq 1$.}
\enq

\begin{claim}\label{Mbound}
$\sum_{i=1}^s M_i \leq 2Lm/\gl^2.$
\end{claim}

We need the following easy observation. 
For $\un{p} = (p_i)_{i=1}^s$ with $p_i \in [0,1)$, let \[f(\un{p}) = \sum_{i=1}^s p_i \prod_{j<i} (1-p_i)^{-1}, \qquad g(\un{p}) = \prod_{j=1}^s (1-p_i)^{-1}.\]
\begin{prop}\label{tech}
$f(\un{p}) \leq g(\un{p})-1$.
\end{prop}
\begin{proof}
We prove the equivalent 
\[\prod_{i=1}^s (1-p_i) + \sum_{i=1}^s p_i\prod_{j\geq i} (1-p_i) \leq 1\]
by induction on $s\geq 1$. 
The base case is obvious, and for the induction step we just observe that
the left hand side is
\[(1-p_1)^2\prod_{i=2}^s (1-p_i) + \sum_{i=2}^s p_i \prod_{j \geq i} (1-p_i) \leq \prod_{i=2}^s (1-p_i) + \sum_{i=2}^s p_i \prod_{j \geq i} (1-p_i) \leq 1.\]
\end{proof}
\begin{proof}[Proof of Claim~\ref{Mbound}]

With sums over $v$'s (in $W$) with $T_v, r_v> i-1$ (cf.\ the discussion following \eqref{xisxis}), we have
\beq{break}
M_i = \E|\sum( \xi_i^v - \xi_{i-1}^v)|
\leq \sum \E|\xi_i^v - \xi_{i-1}^v|= \sum 2\vp(x_{i-1}, v)\xi_{i-1}^v. \\
\enq
Thus, using Proposition~\ref{tech}
for \eqref{useprop}, we have
\begin{align}
    \sum_{i=1}^s M_i &\leq 2\sum_{i=1}^s \sum_v \vp(x_{i-1},v)\xi_{i-1}^v\nonumber\\
    &=2 \sum_v \sum\{\vp(x_{i-1},v)\xi_{i-1}^v:i-1 < T_v\wedge r_v\}\nonumber\\
    &\leq\label{useprop} 2\sum_v (H_v(T_v \wedge r_v)^{-1} - 1)\\
    &\leq 2mL/\gl^2\nonumber.\qedhere
\end{align}
\end{proof}
\begin{proof}[Proof of Theorem~\ref{alltime}]
Lemma~\ref{adaptive}, with \eqref{product} and Claim~\ref{Mbound}, gives
\[\E[e^{\gz Z}] \leq e^{16 L^2 m \gz^2/\gl^4 }
\]
whenever $|\gz| \leq \gl^2/L$. So for any $\gz \in (0, \gl^2/L]$,
\[
\pr(Z > \vt m) = \pr(e^{\gz Z} > e^{\gz \vt m}) \leq  \exp[16 L^2m\gz^2/\gl^4  - \gz \vt m].\]
Since $\xi_s<(L/\gl)m$, 
we may assume $\vt \le L/\gl$ (or the theorem is trivial). 
Setting $\gz = \vt\gl^4/(32L^2)$ to minimize the exponent, we have
\[
\pr(Z > \vt m)\leq \exp[-\vt^2 \gl^4m/(8L)^2 ].
\]
Similarly
\[
\pr(Z < -\vt m) \leq \exp[-\vt^2\gl^4m/(8L)^2 ],
\]
completing the proof.
\end{proof}


\begin{thebibliography}{99}


\bibitem{Aldous}  D.\ Aldous, Lower bounds for covering times for reversible Markov chains and random walks on graphs, \textit{J Theor. Probab.} \textbf{2} (1989), 91–100.

\bibitem{AF} D.\ Aldous and J.\ Fill, \textit{Reversible Markov Chains and Random Walks on Graphs}, Unfinished monograph, available at \href{http://www.stat.berkeley.edu/~aldous/RWG/book.html}{http://www.stat.berkeley.edu/~aldous/RWG/book.html} 2002.



\bibitem{IB} I.\ Benjamini, personal communication.

\bibitem{BGM}  I.\ Benjamini, O.\ Gurel-Gurevich, and B.\ Morris, Linear cover time is exponentially unlikely, \textit{Probab. Theory Relat. Fields} \textbf{155} (2013), 451–461.

\bibitem{MGT} B. Bollob\'as, \textit{Modern Graph Theory}, Springer-Verlag, New York, 1998.

\bibitem{Ding}  J.\ Ding, 
Asymptotics of cover times via Gaussian free fields: Bounded-degree graphs
and general trees, \emph{Ann.\ Probab.} \textbf{42} (2014), 464–496.

\bibitem{DLP}  J.\ Ding, J.\ Lee and Y. Peres,  
Cover times, blanket times, and majorizing measures,
\emph{Ann.\ Math.} \textbf{175} (2012), 1409-1471.

\bibitem{Gillman}  D.\ Gillman, A Chernoff bound for random walks on expander graphs, \textit{SIAM J. Comput.} \textbf{27} (1998), 1203–1220.

\bibitem{aglc} J.\ Kahn, Asymptotically good list-colorings, \textit{J. Combin. Theory Ser. A}
\textbf{73} (1996), 1-59.

\bibitem{MCMT} D.\ Levin, Y.\ Peres, and E.\ Wilmer, \textit{Markov Chains and Mixing Times}, American Mathematical Society, Providence, 2017. With a chapter by James G. Propp and David B. Wilson.

\bibitem{Lovasz}  L.\ Lov\'asz, Random walks on graphs: a survey, \textit{Combinatorics, Paul Erd\H{o}s is eighty} \textbf{2} (1993), 1-46.

\bibitem{KT} S.\ Karlin and H.\ Taylor, \textit{A First Course in Stochastic Processes}, Academic Press, New York, 1975.

\bibitem{Yehudayoff}  A.\ Yehudayoff, Linear cover time for trees is exponentially unlikely, \textit{Chic. J. Theor. Comput. Sci.} \textbf{2012} (2012).

\bibitem{Zhai}  A.\ Zhai, 
Exponential concentration of cover times, 
\emph{Electron. J. Probab.} \textbf{23} (2018), Paper No.\ 32, 22 pp.



\end{thebibliography}
\end{document}